\documentclass[12pt, a4paper]{article}

\usepackage[T1]{fontenc}
\usepackage[latin2]{inputenc}
\usepackage{amsmath}
\usepackage{amssymb}
\newtheorem{thm}{Theorem}

\newtheorem{con}{Conjecture}

\begin{document}

\title{On additive complement of a finite set}
\author{Sándor Z. Kiss \thanks{Institute of Mathematics, Budapest University of Technology and Economics, H-1529 B.O. Box, Hungary; Computer and Automation Research Institute of the Hungarian Academy of Sciences, Budapest H-1111, L\'agym\'anyosi street 11; kisspest@cs.elte.hu. This author was supported by the OTKA Grant No. K77476 and No. NK105645.}, Eszter Rozgonyi \thanks{Institute of Mathematics, Budapest University of Technology and Economics, H-1529 B.O. Box, Hungary, reszti@math.bme.hu. The work reported in the paper has been developed in the framework of the project "Talent care and cultivation in the scientific workshops of BME" project. This project is supported by the grant TÁMOP - 4.2.2.B-10/1--2010-0009.}, Csaba Sándor \thanks{Institute of Mathematics, Budapest University of Technology and Economics, H-1529 B.O. Box, Hungary, csandor@math.bme.hu. This author was supported by the OTKA Grant No. K81658.}
}

\date{\today}

\maketitle

\begin{abstract}
\noindent We say the sets of nonnegative integers $\mathcal{A}$ and
$\mathcal{B}$ are additive complements if their sum contains all sufficiently large integers. In this paper we prove a conjecture of Chen and Fang about additive complement of a finite set.
\end{abstract}

\textit{2000 AMS \ Mathematics subject classification number}:
   primary: 11B13, secondary: 11P99.


\section{Introduction}

Let $\mathbb{N}$ denote the set of positive integers and let $\mathcal{A}
\subseteq \mathbb{N}$ and $\mathcal{B} \subseteq \mathbb{N}$ be finite or
infinite sets. Let $R_{\mathcal{A}+\mathcal{B}}(n)$ denote the number of solutions of the equation
\[
a + b = n, \hspace*{3mm} a \in \mathcal{A},\hspace*{3mm} b \in \mathcal{B}.
\]
\noindent We put
\[
A(n) = \sum_{\overset{a \le n}{a \in \mathcal{A}}}1 \hspace*{3mm}
and \hspace*{3mm} B(n) = \sum_{\overset{b \le n}{b \in \mathcal{B}}}1
\]
\noindent respectively. We say a set $\mathcal{B} \subseteq \mathbb{N}$ is an additive
complement of the set $\mathcal{A} \subseteq \mathbb{N}$ if every sufficiently large
$n \in \mathbb{N}$ can be represented in the form $a + b = n$, $a \in
\mathcal{A}$, $b \in \mathcal{B}$, i.e., $R_{\mathcal{A}+\mathcal{B}}(n) \ge
1$ for $n \ge n_{0}$. Additive complement is an important concept in additive
number theory, in the past few decades it was studied by many authors [4], [6],
[8], [9]. In [8] Sárközy and Szemerédi proved a conjecture of Danzer [4],
namely they proved that for infinite additive complements $\mathcal{A}$ and
$\mathcal{B}$ if
\[
\limsup_{x \rightarrow +\infty}\frac{A(x)B(x)}{x} \le 1,
\]
\noindent then
\[
\liminf_{x \rightarrow +\infty}(A(x)B(x) - x) = +\infty.
\]
\noindent In [1] Chen and Fang improved this result and they proved that if  \[
\limsup_{x \rightarrow +\infty}\frac{A(x)B(x)}{x} > 2, \hspace*{3mm}
or \hspace*{3mm} \limsup_{x \rightarrow
   +\infty}\frac{A(x)B(x)}{x} < \frac{5}{4},
\]
\noindent then
\[
\lim_{x \rightarrow +\infty}(A(x)B(x) - x) = +\infty.
\]
\noindent In the other direction they proved in [2] that for any integer $a \ge
2$, there exist two infinite additive complements $\mathcal{A}$ and
$\mathcal{B}$ such that
\[
\limsup_{x \rightarrow +\infty}\frac{A(x)B(x)}{x} = \frac{2a+2}{a+2},
\]
\noindent but there exist infinitely many positive integers $x$ such that
$A(x)B(x) - x = 1$. In [3] they studied the case when $\mathcal{A}$ is a finite
set. In this case the situation is different from the infinite case. Chen and
Fang proved that for any two additive complements $\mathcal{A}$ and
$\mathcal{B}$ with $|\mathcal{A}| < +\infty$ or $|\mathcal{B}| < +\infty$, if

\[
\limsup_{x \rightarrow +\infty}\frac{A(x)B(x)}{x} > 1,
\]
\noindent then
\[
\lim_{x \rightarrow +\infty}(A(x)B(x) - x) = +\infty.
\]
\noindent They also proved that if
\[
\mathcal{A} = \{a + im^{s} + k_{i}m^{s+1}:i = 0,...,m - 1\},
\]
\noindent where $|\mathcal{A}| = m$, $a$, $s \ge 0$ and $k_{i}$ are integers,
then there exists an additive complement $\mathcal{B}$ of $\mathcal{A}$ such
that $A(x)B(x) - x = O(1)$. In the special case $|\mathcal{A}| = 3$ they
proved that if
$\mathcal{A}$ is not of the form $\{a + i3^{s} + k_{i}3^{s+1}:i = 0, 1, 2\}$,
where $a$, $s \ge 0$ and $k_{i}$ are integers, then for any additive complement
$\mathcal{B}$ of $\mathcal{A}$,
\[
\lim_{x \rightarrow +\infty}(A(x)B(x) - x) = +\infty
\]
\noindent holds. Chen and Fang posed the following conjecture (Conjecture 1.5. in [3]):
\begin{con}
If the set of nonnegative integers $\mathcal{A}$ is not of the form
\[
\mathcal{A} = \{a + im^{s} + k_{i}m^{s+1}:i = 0,...,m - 1\},
\]
\noindent where $a, m > 0$, $s \ge 0$ and $k_{i}$ are integers, then, for any
additive complement $\mathcal{B}$ of $\mathcal{A}$, we have
\[
\lim_{x \rightarrow +\infty}(A(x)B(x) - x) = +\infty.
\]
\end{con}
\noindent In this paper we prove this conjecture, when the number of elements
of the set $\mathcal{A}$ is prime:

\begin{thm}
Let $p$ be a positive prime and $\mathcal{A}$ is a set of nonnegative
integers with $|\mathcal{A}| = p$. If $\mathcal{A}$ is not of the form
\begin{equation}
\mathcal{A} = \{a + ip^{s} + k_{i}p^{s+1}:i = 0,...,p - 1\},
\end{equation}
\noindent where $a > 0$, $s \ge 0$ and $k_{i}$ are integers, then, for any
additive complement $\mathcal{B}$ of $\mathcal{A}$, we have
\begin{equation}
\lim_{x \rightarrow +\infty}(A(x)B(x) - x) = +\infty.
\end{equation}
\end{thm}

\noindent In the case when the number of elements of $\mathcal{A}$ is a
composite number, we disprove the Conjecture 1.5. in [3]:

\begin{thm}
For any composite number $n > 0$, there exists a set $\mathcal{A}$ and
a set $\mathcal{B}$ such that $|\mathcal{A}| = n$, $\mathcal{B}$ is an
additive complement of $\mathcal{A}$ and $\mathcal{A}$ is not of the form
\[
\mathcal{A} = \{a + in^{s} + k_{i}n^{s+1}:i = 0,...,n - 1\},
\]
\noindent where $s \ge 0$, $a > 0$, and $k_{i}$ are integers, and
\[
A(x)B(x) - x = O(1).
\]
\end{thm}
\noindent In the next section we give a short survey about the algebraic
concepts which play a crucial role in the proof of Theorem 1.

\section{Preliminaries}
In our proof we are working with cyclotomic polynomials. Both the definition
and the most important properties of these polynomials are
well-known. Interested reader can find these in [5, p. 63-66]. We denote the
degree of a polynomial $f$ by $\deg f$. Let $\theta$ be an algebraic
number. We say the monic polynomial $f$ is the minimal polynomial of $\theta$
if $f$ is the least degree such that $f(\theta) = 0$. It is well-known that if $f$ is the minimal
polynomial of $\theta$, and $g$ is a polynomial such that $g(\theta) = 0$,
then $f|g$. A $\mu$ complex number is called primitive $n$th root of unity if
$\mu$ is the root of the polynomial $x^{n} - 1$ but not of $x^{m} - 1$ for any
$m < n$. The cyclotomic polynomial of order $n$ is defined by
\[
\Phi_{n}(z) = \prod_{\zeta}(z - \zeta),
\]
\noindent where $\zeta$ runs over all the primitive $n$th root of
unity. This is a monic irreducible polynomial with degree $\varphi(n)$, and
$\Phi_{n}(z)$ has integer coefficients. It is well-known that $\Phi_{n}(z)$ is
the minimal polynomial of $\zeta$ and

\begin{equation}
1 + z + z^{2}+ \dots{} + z^{n-1} = \prod_{\overset{l|n}{l > 1}}\Phi_{l}(z).
\end{equation}
\noindent It is easy to see that

\begin{equation}
\Phi_{p^{s+1}}(z) = 1 + z^{p^{s}} + z^{2p^{s}} + \dots{} + z^{(p-1)p^{s}}
\end{equation}

\section{Proof of Theorem 1}

We will prove that if there exists an additive complement
$\mathcal{B}$ of $\mathcal{A}$, $|\mathcal{A}| = p$ such that
\[
\liminf_{x \rightarrow +\infty}(A(x)B(x) - x) < +\infty,
\]
\noindent then $\mathcal{A}$ is the form (1). Let us suppose that $R_{\mathcal{A}+\mathcal{B}}(n)\geq 1$ for $n\geq n_0$. First we prove that there exists
an integer $n_{1}$ such that $R_{\mathcal{A}+\mathcal{B}}(n) = 1$ for $n \ge
n_{1}$. We argue as Sárközy and Szemerédi in [9, p.238]. As $\mathcal{B}$ is an
additive complement of $\mathcal{A}$, it follows that
\[
+\infty > C = \liminf_{x \rightarrow +\infty}(A(x)B(x) - x) =
\liminf_{x \rightarrow +\infty}\Bigg(\Big(\sum_{\overset{a \in
   \mathcal{A}}{a \le x}}1\Big)\Big(\sum_{\overset{b \in \mathcal{B}}{b \le
   x}}1\Big) - x\Bigg) \ge
\]
\[
\ge \liminf_{x \rightarrow +\infty}\Bigg(\Big(\sum_{\overset{a \in
   \mathcal{A}, b \in \mathcal{B}}{a + b \le x}}1\Big) - x\Bigg) =
\liminf_{x \rightarrow
   +\infty}\Big(\sum_{n=0}^{x}R_{\mathcal{A}+\mathcal{B}}(n) - x\Big) \ge
\]
\[
\ge \liminf_{x \rightarrow
   +\infty}\Big(\sum_{n=n_{0}+1}^{x}R_{\mathcal{A}+\mathcal{B}}(n) - x\Big)
\ge \liminf_{x \rightarrow +\infty}\Bigg([x] - n_{0} +
\sum_{\overset{n_{0} < n \le x}{R_{\mathcal{A}+\mathcal{B}}(n) > 1}}1 - x\Bigg) \ge
\]
\[
\ge \liminf_{x \rightarrow +\infty}\Bigg(
\sum_{\overset{n_{0} < n \le x}{R_{\mathcal{A}+\mathcal{B}}(n) > 1}}1\Bigg) -
(n_{0} + 1),
\]
\noindent thus we have
\[
\liminf_{x \rightarrow +\infty}\Bigg(
\sum_{\overset{n_{0} < n \le x}{R_{\mathcal{A}+\mathcal{B}}(n) > 1}}1\Bigg) < C + n_{0} + 1,
\]
\noindent where $C$ is a positive constant. As $\mathcal{B}$ is an additive complement of $\mathcal{A}$, it
follows that there exists an integer $n_{1}$ such that
\begin{equation}
R_{\mathcal{A}+\mathcal{B}}(n) = 1 \hspace*{3mm} for \hspace*{3mm} n \ge n_{1}.
\end{equation}
In the next step we prove that $\mathcal{A}$ is the form (1). Let $z = re^{2i\pi\alpha} = re(\alpha)$, where $r < 1$. Let the generating
functions of the sets $\mathcal{A}$ and $\mathcal{B}$ be $f_{\mathcal{A}}(z)
= \sum_{a \in \mathcal{A}}z^{a}$ and $f_{\mathcal{B}}(z) = \sum_{b \in
  \mathcal{B}}z^{b}$ respectively. (By $r < 1$ these infinite series and all
the other infinite series of the proof are absolutely convergent.) In view of (5) we have
\[
f_{\mathcal{A}}(z) \cdot f_{\mathcal{B}}(z) = \Big(\sum_{a \in \mathcal{A}}z^{a}\Big)\Big(\sum_{b \in
  \mathcal{B}}z^{b}\Big) =
\sum_{n=0}^{+\infty}R_{\mathcal{A}+\mathcal{B}}(n)z^{n} =
\]
\[
= \sum_{n=0}^{n_{1}- 1}R_{\mathcal{A}+\mathcal{B}}(n)z^{n} +
\sum_{n=n_{1}}^{+\infty}R_{\mathcal{A}+\mathcal{B}}(n)z^{n} = p_{1}(z) + \frac{z^{n_{1}}}{1 - z},
\]
\noindent where $p_{1}(z)$ is a polynomial of $z$. Thus we have
\begin{equation}
(1 - z)f_{\mathcal{A}}(z) \cdot f_{\mathcal{B}}(z) = (1 - z)p_{1}(z) + z^{n_{1}}.
\end{equation}
\noindent In next step we prove that $f_{\mathcal{B}}(z)$ can be written in
the form
\begin{equation}
f_{\mathcal{B}}(z) = F_{\mathcal{B}}(z) + \frac{T(z)}{1 - z^{M}},
\end{equation}
\noindent where $M$ is a positive integer, $F_{\mathcal{B}}(z)$ and $T(z)$ are
polynomials. We argue as Nathanson in [7, p.18-19]. Let $(1 - z)f_{\mathcal{A}}(z) = \sum_{n=K}^{N}a_{n}z^{n}$, where
$a_{N} \ne 0$ and $a_{K} \ne 0$, and let $f_{\mathcal{B}}(z) =
\sum_{n=0}^{\infty}e_{n}z^{n}$, where $e_{n} \in \{0, 1\}$.  Then we have
\[
(1 - z)f_{\mathcal{A}}(z) \cdot f_{\mathcal{B}}(z) = \sum_{n=0}^{\infty}c_{n}z^{n},
\]
\noindent where $c_{n} = 0$ from a certain point on. It is clear that if $n$
is large enough, then $c_{n} = e_{n-K}a_{K} + e_{n-K-1}a_{K+1} + \dots{} +
e_{n-N}a_{N} = 0$. This shows that the coefficients of the power series
$f_{\mathcal{B}}(z)$ satisfies a linear recurrence relation from a certain
point on. These coefficients are either $0$ or $1$ from a certain point on. It
is easy to see that a sequence defined by a linear recurrence relation on a
finite set must be eventually periodic, which proves (7). 

\noindent It follows from (6) and (7) that

\[
f_{\mathcal{A}}(z) \cdot \Bigg(F_{\mathcal{B}}(z) + \frac{T(z)}{1 -
 z^{M}}\Bigg) = p_{1}(z) + \frac{z^{n_{1}}}{1 - z},
\]
\noindent hence for every $z\in \mathbb{C}$
\begin{equation}
(1 - z^{M})f_{\mathcal{A}}(z)F_{\mathcal{B}}(z) + f_{\mathcal{A}}(z)T(z) = (1
 - z^{M})p_{1}(z) + (1 + z + z^{2}+ \dots{} + z^{M-1})z^{n_{1}}.
\end{equation}
By putting $z = 1$, we obtain that
\begin{equation}
f_{\mathcal{A}}(1)T(1) = M.
\end{equation}
As $f_{\mathcal{A}}(1) = |\mathcal{A}| = p$, it follows from (9) that
$p|M$. Define $k$ by $p^{k}|M$ but $p^{k+1} \nmid M$. It follows from (  that
\[
(1 + z + z^{2} + \dots{} + z^{M-1})|f_{\mathcal{A}}(z)T(z).
\]
It follows from (3) that for any $1 \le t \le k$ we have
\[
\Phi_{p^{t}}(z)|f_{\mathcal{A}}(z)T(z).
\]
Assume that for any $1 \le t \le k$ we have $\Phi_{p^{t}}(z)|T(z)$. Then
\[
T(z) = \left( \prod_{t = 1}^{k}\Phi_{p^{t}}(z)\right) \cdot q(z),
\]
\noindent where $q(z)$ is a polynomial with integer coefficients. By putting $z = 1$ we obtain that $T(1) = p^{k}q(1)$, hence $M = f_{\mathcal{A}}(1)T(1) =
p^{k+1}q(1)$ which contradicts the definition of $k$. It follows that
there exists an integer $0
\le s \le k - 1$ such that $\Phi_{p^{s+1}}(z)|f_{\mathcal{A}}(z)$, thus
$f_{\mathcal{A}}(z) = \Phi_{p^{s+1}}(z) \cdot a(z)$, where $a(z)$ is a polynomial. As $\mathcal{A} =
\{a_{1}, \dots{}, a_{p}\}$, we have $f_{\mathcal{A}}(z)
= \sum_{i=1}^{p}z^{a_{i}}$. Let $\omega$ be the following $p^{s+1}$th root of unity,
\[
\omega = exp\Bigg(\frac{2\pi}{p^{s+1}}i\Bigg).
\]
\noindent It follows that
$f_{\mathcal{A}}(\omega) = 0$, thus we have $\sum_{i=1}^{p}\omega^{a_{i}} =
0$. Let $a_{i} = l_{i}p^{s+1} + r_{i}$, where $0 \le r_{i} < p^{s+1}$. Without
loss of generality we may assume that
\begin{equation}
0 \le r_{1} \le r_{2} \le \dots{} \le r_{p} < p^{s+1}.
\end{equation}
Define $r_{p+1} = p^{s+1}
+ r_{1}$. Since $\sum_{i=1}^{p}(r_{i+1} - r_{i}) = r_{p+1} - r_{1} = p^{s+1}$
then it follows that there exists a $j$ with $1 \le j \le p$ such that
\begin{equation}
r_{j+1} - r_{j} \ge p^{s}.
\end{equation}
In the next step we prove that this implies
\begin{equation}
a_{i} - r_{j+1} = n_{i}p^{s+1} + t_{i},
\end{equation}
\noindent where $1 \le i \le p$ and  $0 \le t_{i} \le p^{s+1} - p^{s}$ holds.
Assume that $i \le j$. By the definition of $a_{i}$ we have
$a_{i} - r_{j+1} = l_{i}p^{s+1} + r_{i} - r_{j+1}$. It follows from (10) and
(11) that $r_{j+1} - r_{i} \le r_{j+1} < p^{s+1}$ and
$-p^{s+1} < r_{i} - r_{j+1} \le r_{j} - r_{j+1} \le -p^{s}$. Thus we have $0
\le r_{i} - r_{j+1}+p^{s+1} \le p^{s+1} - p^{s}$, which implies (12).
In the second case assume that $i \ge j + 2$. It is clear from
(10) that $r_{i} - r_{j+1} > 0$. By the definition of $a_{i}$
and (10), (11) we have $$a_{i} - r_{j+1} = l_{i}p^{s+1} + r_{i} -
r_{j+1} <l_{i}p^{s+1}+p^{s+1}-r_{j+1}\le l_{i}p^{s+1} + p^{s+1} - p^{s},$$ which implies (12).
It follows that there exists an integer $a$ such that $a_{i} = a +
n_{i}p^{s+1} + t_{i}$, where $n_{i}$ is an integer and
\begin{equation}
0 \le t_{i} \le p^{s+1} - p^{s}.
\end{equation}
As $f_{\mathcal{A}}(\omega) = 0$, and the definition of $\omega$ we obtain that
\[
\sum_{i=1}^{p}\omega^{a_{i}} =
\sum_{i=1}^{p}\omega^{a + n_{i}p^{s+1} + t_{i}} = \sum_{i=1}^{p}\omega^{a +
 t_{i}} = 0.
\]
\noindent Let $h(z) = \sum_{i=1}^{p}z^{t_{i}}$. Thus we obtain that $h(\omega) =
0$. As $\Phi_{p^{s+1}}(z)$ is a minimal polynomial of $\omega$ we have
$\Phi_{p^{s+1}}(z)|h(z)$. It follows from (13) that
$deg\Big(\sum_{i=1}^{p}z^{t_{i}}\Big) \le p^{s+1} - p^{s} = \varphi(p^{s+1}) =
deg\Big(\Phi_{p^{s+1}}(z)\Big)$. Therefore, by using (4) we have
$\sum_{i=1}^{p}z^{t_{i}} = \Phi_{p^{s+1}}(z) = 1 + z^{p^{s}} + z^{2p^{s}} + \dots{} + z^{(p-1)p^{s}}$ and
then we have $\{t_{1}, \dots{}, t_{p}\} = \{0, p^{s}, 2p^{s}, \dots{},
(p-1)p^{s}\}$. It follows that there exist integers $a > 0$ and $k_{i}$, such
that $\mathcal{A} = \{a + ip^{s} + k_{i}p^{s+1}\}$, as desired.

\section{Proof of Theorem 2}

Let $n = d_{1}d_{2}$, $d_{1}$, $d_{2} > 1$ be integers, and consider the following two sets:
\[
\mathcal{A} = \{u + v \cdot d_{1}d_{2}: 0 \le u \le d_{1} - 1, 0 \le v \le
d_{2} - 1\},
\]
\[
\mathcal{B} = \{kd_{1}d_{2}^{2} + w \cdot d_{1}: k \in \mathbb{N}, 0 \le w \le
d_{2} - 1\}.
\]
\noindent It is easy to see that $|\mathcal{A}| = d_{1}d_{2}$. It is clear
that $A(x) = d_{1}d_{2}$ if $x$ is large enough and $B(x)=
\frac{x}{d_{1}d_{2}} + O(1)$, which
 implies $A(x)B(x)-x=O(1)$. Let $m$ be a fixed positive integer. It is clear that any positive integer $m$ can be written uniquely in the form  
\[
m = kd_{1}d_{2}^{2} + ud_{1} + ld_{1}d_{2} + v,
\]
\noindent where $k$ is a nonnegative integer, $0 \le u,l \le d_{2}$, $0 \le v
\le d_{1}$. Hence $\mathcal{B}$ is an additive complement of
$\mathcal{A}$. In the next step we prove that the set $\mathcal{A}$ is not of
the form (1). Assume that $\mathcal{A}$ is the form (1). It is clear that the
difference of any two elements from $\mathcal{A}$ divisible by $n^{s}$. As
$\mathcal{A}$ also contains consecutive integers we have $n^{s}|1$, which
implies $s = 0$. Thus $\mathcal{A} = \{a + i + k_{i}n: i = 0, \dots{}, n -
1\}$, that is $\mathcal{A}$ is a complete residue system modulo $n$, which is a
contradiction.

\bigskip
\bigskip

\textsc{Acknowledgement}: The authors would like to thank András Bíró for the
valuable discussions.

\end{document}